\documentclass[10pt,twocolumn,twoside]{IEEEtran}

\usepackage{graphicx}
\usepackage{caption}
\usepackage{subfig}
\usepackage{cite}
\usepackage{stmaryrd}
\usepackage{amsfonts,latexsym,amssymb}
\usepackage[cmex10]{amsmath}
\usepackage{algorithm}
\usepackage{algorithmic}
\usepackage{url}

\input{zgx.math}

\newtheorem{theorem}{{\bf{Proposition}}}
\newtheorem{lemma}{{\bf{Lemma}}}
\newtheorem{corollary}{{\bf Corollary}}

\begin{document}
\title{Accelerated Canonical Polyadic Decomposition by Using Mode Reduction}
\author{Guoxu~Zhou \IEEEmembership{Member,~IEEE}, Andrzej~Cichocki ~\IEEEmembership{Fellow,~IEEE}, and Shengli~Xie \IEEEmembership{Senior Member,~IEEE,}

\thanks{Manuscript received December 7, 2012; revised April 1, 2013; accepted June 22, 2013. This work was supported in part by National Natural Science Foundation of China under Grant 61103122, Grant 61273192, Grant U1201253, Grant 61202155, and the Guangdong Natural Science Foundation under Grant S2011040005724.}
\thanks{Guoxu Zhou is with the Laboratory for Advanced Brain Signal Processing, RIKEN, Brain Science Institute, Wako-shi, Saitama 3510198, Japan. He is also with the Faculty of Automation, Guangdong University of Technology, Guangzhou, 510641, China. E-mail: zhouguoxu@ieee.org.}%
\thanks{Andrzej Cichocki is with the RIKEN BSI, Japan and Systems Research Institute, Warsaw Poland. E-mail: cia@brain.riken.jp.}%
\thanks{Shengli Xie is with the Faculty of Automation, Guangdong University of Technology, Guangzhou 510006, China. E-mail: eeoshlxie@scut.edu.cn.}
}

\markboth{IEEE TRANSACTIONS ON NEURAL NETWORKS AND LEARNING SYSTEMS}%
{ZHOU \MakeLowercase{\textit{et al.}}: ACCELERATED CANONICAL POLYADIC DECOMPOSITION BY USING MODE REDUCTION}

\maketitle

\begin{abstract}
  Canonical Polyadic (or CANDECOMP/PARAFAC, CP) decomposition (CPD) is widely applied to $N$th-order ($N\ge3$) tensor analysis. Existing CPD methods mainly use alternating least squares iterations and hence need to unfold tensors to each of their $N$ modes frequently, which is one major performance bottleneck for large-scale data, especially when the order $N$ is large. To overcome this problem, in this paper we proposed a new CPD method in which the CPD of a high order tensor (i.e. $N>3$) is realized by applying CPD to a mode reduced one (typically, 3rd-order tensor) followed by a Khatri-Rao product projection procedure. This way is not only quite efficient as frequently unfolding to $N$ modes is avoided, but also promising to conquer the bottleneck problem caused by high collinearity of components. We showed that, under mild conditions, any $N$th-order CPD can be converted to an equivalent 3rd-order one but without destroying essential uniqueness, and theoretically they simply give consistent results. Besides, once the CPD of any unfolded lower-order tensor is essentially unique, it is also true for the CPD of the original higher order tensor. Error bounds of truncated CPD were also analyzed in presence of noise. Simulations showed that, compared with state-of-the-art CPD methods, the proposed method is more efficient and is able to escape from local solutions more easily.
\end{abstract}


% \fi
\IEEEpeerreviewmaketitle

\section{Introduction}
\label{sec:Intro}
\IEEEPARstart{H}{igher-order} tensors (multi-way arrays) have gained increasing importance as they are often more natural representations of multi-dimensional data than matrices in many practical applications \cite{NMF-book,Kolda09tensordecompositions}. As one of the most fundamental problem in tensor data analysis, tensor decomposition attempts to find informative representations (e.g. dense/sparse, low-rank representation) of multi-dimensional tensor data. Tensor decomposition is very attractive and versatile because it takes  into account such as spatial, temporal and spectral information, and provides links among the various extracted  factors or latent variables with desired physical or physiological meaning and interpretation \cite{NMF-book,Kolda09tensordecompositions,tnn_tensorDisc,TNN2009UMDA,SpTDA_FaceVeri,TNNTensorDist}.

As one of the most important tensor decomposition models, Canonical Polyadic
(CP), also named as CANDECOMP/PARAFAC decomposition \cite{PARAFAC1970Carroll,PARAFAC1970Harsman}, has been extensively studied in the last four decades and found many practical applications \cite{NMF-book}, for example, in underdetermined blind source separation \cite{TNN-TFUBSS}. One major advantage of CPD is that it is essentially unique under mild conditions (see Lemma \ref{lemma:Uniqueness} in Section \ref{sec:Uniqueness}), which makes it very useful in the case where only very limited or even no \apriori knowledge is available on factors. In the CP model, the mode-$n$ matricization of a tensor is just a product of its mode-$n$ factor matrix with another matrix formed by the Khatri-Rao product of all the remaining matrices in a specific order. This feature has led to the widely adopted alternating least squares (ALS) methods to solve CPD, e.g. see \cite{KoldaTensorToolbox}. Unfortunately, these methods require to unfold the tensor to its $N$ modes frequently, which is one major performance bottleneck of CPD algorithms.

In the CPD of an $N$th-order tensor there are a total of $N$ factor matrices to be estimated. Surprisingly, our recent results showed that once at least one factor with full column rank has been correctly estimated, all the other factors can be computed uniquely and efficiently by using a series of singular value decompositions (SVD) of rank-1 matrices\footnote{A matrix \mat{Y} is rank-1 if and only if $\mat{Y}=\mat{uv}^T$, where \mat{u} and \mat{v} are two nonzero vectors.}, no matter whether they contain some collinear components (columns) \cite{SPL-smbss}. This motivated us to perform blind source separation (BSS) on one single mode first and then use efficient rank-1 approximation methods to recover the  other factors, which has led to the CP-SMBSS method for CP decompositions. The CP-SMBSS is very useful if some \apriori knowledge about component diversities, such as independence \cite{Common1994,MaartenDeVos2008}, nonnegativity \cite{Lee1999,TNN-MVCNMF,Yang_TNN2012}, sparsity \cite{TNN-sMixEst}, etc, in at least one mode is available, which allows us to recover them from their linear mixtures via standard BSS algorithms. In this paper, however, we focus on the case where such \apriori information is completely unavailable.

The following notations will be adopted. Bold capitals (e.g., $\mat{A}$) and bold lowercase letters (e.g., $\mat{y}$) denote matrices and vectors, respectively. Calligraphic bold capitals, e.g. \tensor{Y}, denote tensors. Mode-$n$ matricization (unfolding, flattening) of a tensor  $\tensor{Y} \in \Real^{I_{1} \times I_{2} \times \cdots \times I_{N}}$ is denoted as $\tenmat{Y} \in \Real^{I_{n} \times \prod_{p\neq{n}}{I_p}}$, which consists of arranging all possible mode-$n$ tubes (vectors) as the columns of it \cite{Kolda09tensordecompositions}. The Frobenius norm of a tensor is denoted by $\frob{\tensor{Y}}=(\sum_{i_1i_2\cdots i_N}y_{i_1i_2\cdots i_N}^2)^{\frac{1}{2}}$.

The mode-$n$  product of a tensor
$\tensor{G} \in \Real^{J_{1} \times J_{2} \times \cdots \times J_{N}}$ and
a matrix $\mat{A} \in \Real^{I \times J_n}$ yields  a tensor denoted by
$\tensor{Y}=\tensor{G} \ttm{\mat{A}} \in \Real^{J_1 \times \cdots \times J_{n-1} \times I \times J_{n+1} \times \cdots \times J_N}$, with elements $y_{j_1,j_2,\ldots,j_{n-1},i,j_{n+1},\ldots, j_{N}} =\sum_{j_n=1}^{J_n} (g_{j_1,j_2,\ldots,j_N})(a_{i,j_n})$.

We use $\krp$ and $\hdp$ to denote the Khatri-Rao product (column-wise Kronecker product) and Hadamard product of matrices, respectively. Given a set of matrices $\matn[k]{A}\in\Real^{I_k\times J}$ with  $k=k_1,k_1+1,k_1+2,\ldots,k_2$, $\bigkrp\nolimits_{k=k_1}^{k_2}\matn[k]{A}=\matn[k_2]{A}\krp\matn[k_2-1]{A}\cdots\krp\matn[k_1+1]{A}\krp\matn[k_1]{A}$ .
Readers are referred to \cite{Kolda09tensordecompositions,NMF-book} for detailed tensor notations and operations.

\section{Model Reduction Of Tensors In CP Decompositions}
\label{sec:MR}
\subsection{CP Decomposition of Tensors}
\label{subsec:CPD}
\begin{figure}[!t]
\centerline{
\includegraphics[width=.49\textwidth]{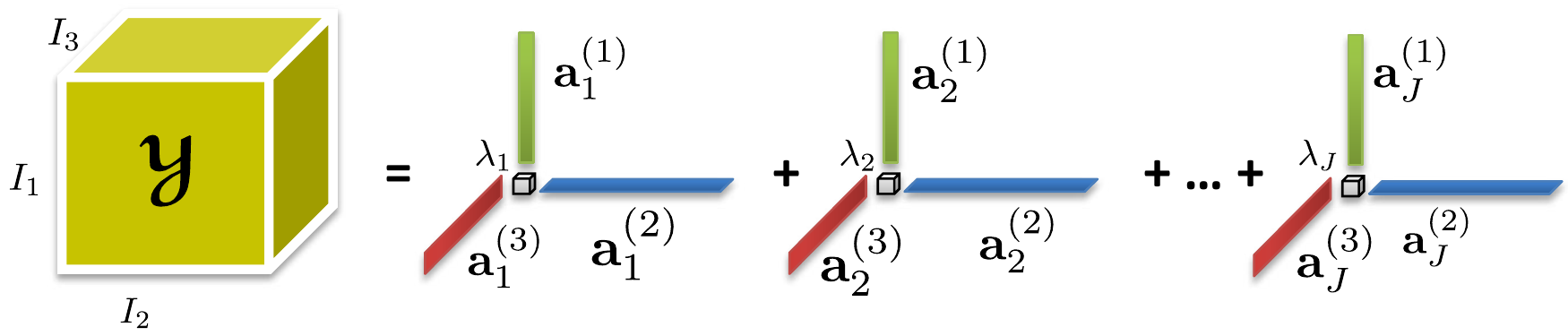}
}
\caption{Illustration of CP decompositions of a 3rd-order tensor $\tensor{Y}\in\Real^{I_1\times I_2\times I_3}$ (ignored the noise), where the factors $\matn{A}=[\matn{a}_1 \; \matn{a}_2 \; \cdots \; \matn{a}_J]\in\Real^{I_n\times J}$ contain the latent components $\matn{a}_j$ as their columns, $n=1,2,3$. }
\label{figcptensor}
\end{figure}

 CP decomposition (or factorization) of a tensor \tensor{Y}$\in\Real^{I_1\times I_2\cdots\times I_N}$ can be  formulated as
\begin{equation}
\label{eqGeneralCPModel}
\begin{split}
 \tensor{Y} &= \sum_{j=1}^J  \lambda_j \; \matn[1]{a}_j \circ  \matn[2]{a}_j \cdots \circ \matn[N]{a}_j +\tensor{E},\\
           %%        &=\compactcp{A}+\tensor{E},
\end{split}
\end{equation}
where component (or factor, mode) matrices $\matn{A} = [ \matn{a}_1,\matn{a}_2, \cdots, \matn{a}_J] \in \Real^{I_n \times J}$, $n\in\Set{N}=\set{1,2,\cdots,N}$, consist of unknown latent components $\matn{a}_j$ (e.g., latent source signals) that need to be estimated, $\outerp$ denotes the outer product\footnote{The outer product of two vectors
$\mat{a} \in \Real^I, \; \mat{b} \in \Real^T$ builds up a rank-one matrix
$\mat{Y}=\mat{a} \outerp \mat{b} = \mat{a} \mat{b}^T \in \Real^{I \times T}$
and the outer product of  three vectors: $\mat{a} \in \Real^I, \; \mat{b} \in \Real^T, \; \mat{c} \in \Real^Q$ builds up a 3rd-order rank-one tensor:
$\tensor{Y}= \mat{a} \outerp \mat{b} \outerp \mat{c} \in \Real^{I \times T  \times Q}$,
with entries defined as $y_{itq} =a_i b_t c_q$.},
 and \tensor{E} denotes the tensor of error or residual terms. Note that \tensor{E} can be zero if we  increase the value of $J$ arbitrarily. In fact, any $N$th-order tensor \tensor{Y} with finite dimension can always be exactly represented by using a CP model with some finite $J$, since $\tensor{Y}=\sum_{i_1,i_2,\ldots,i_N}y_{i_1 i_2 \cdots i_N}\matn[1]{e}_{i_1}\outerp \matn[2]{e}_{i_2}\cdots\outerp \matn[N]{e}_{i_N}$, where the vector $\matn[n]{e}_{i_n}\in\Real^{I_n\times 1}$ has all zero entries except its $i_n$th entry equal to 1.  In CPD the minimum value of $J$ is of particular interest and is called the rank of tensor \tensor{Y}. See \figurename  \ref{figcptensor} for the illustration of the CPD of a 3rd-order tensor. From (\ref{eqGeneralCPModel}), the multiway tensor is represented as a linear combination of outer products of vectors (i.e., rank one tensors), which can be regarded as a generalization of  matrix SVD to the tensor case \cite{Kolda09tensordecompositions}. As the scalar factors $\lambda_j$ can be absorbed into one factor matrix, e.g. \matn[N]{A} by letting $\matn[N]{a}_j=\lambda_j\matn[N]{a}_j, \forall j$, we also use $\tensor{Y}=\compactcp{A}$ as a shorthand notation of (\ref{eqGeneralCPModel}). For this reason hereafter we assume that $\frob[2]{\matn{a}_j}=1$ for all $j=1,2,\ldots,J$ and $n\neq N$.

To solve CPD problems, alternating least squares (ALS) methods are widely employed. Consider the mode-$n$ matricization of \tensor{Y}:
\begin{equation}
\label{CPmatn}
  \tenmat{Y}=\matn{A}{\matn{B}}{}^T, \;(n\in\Set{N}),
\end{equation}
where
\begin{equation}
  \label{CPBn}
  \matn{B}=\bigkrp\nolimits_{p\neq n}\matn[p]{A}\in\Real^{(\prod_{p\neq n}{I_p})\times J}.
\end{equation}
In standard ALS based algorithms, factor matrices \matn{A} are updated using $\matn{A}\from \tenmat{Y}[{\matn{B}{}^{T}]^\dagger}$ alternatively for $n=1,2,\ldots,N$, where \pinv{} denotes the Moore-Penrose pseudo inverse of a matrix. As the matrix \matn{B} is often quite huge, some tricks were proposed to simplify the computation of $\tenmat{Y}[{\matn{B}{}^{T}]^\dagger}$, for example, see \cite{Kolda09tensordecompositions}.

One of the most attractive property of CPD is that it is essentially unique under mild conditions \cite{Kruskal1977,Sidiropoulos2000} (see also Lemma \ref{lemma:Uniqueness} in section \ref{sec:Uniqueness}), which means that for another CPD of \tensor{Y} such that $\tensor{Y}=\tenfactors{\matn[1]{\widehat{A}},\matn[2]{\widehat{A}},\ldots,\matn[N]{\widehat{A}}}$, there must hold that
\begin{equation}
  \matn{\widehat{A}}=\matn{A}\mats{P}\mats{D}, \forall n\in\Set{N},
\end{equation}
where \mats{P} and \mats{D} are any permutation matrices and nonsingular diagonal matrices associated with mode-$n$, respectively.

\subsection{Mode Reduction of Tensors}
\label{subsec:MR}
From the above analysis, when the number of modes $N$ is large, ALS methods often suffer from very slow convergence speed as they need to unfold  tensors with respect to each of the $N$ modes frequently. A number of authors have made efforts to improve the efficiency of CPD algorithms, e.g. see  \cite{PComon2009,ELS::LDL2011}.
In this paper we consider a new way to conquer this problem, i.e., reducing the number of modes to accelerate the convergence of CPD algorithms. As preliminary we need the following tensor operations.

{\bf Tensor transpose}\cite{TensorUnfoldings}. Given an $I_1\times I_2 \cdots \times I_N$ tensor \tensor{Y}=\compactcp{A}  defined in (\ref{eqGeneralCPModel}), the transpose of \tensor{Y} is a tensor of ${I_{p_1}\times I_{p_2} \cdots \times I_{p_N}}$ obtained by exchanging the roles of \matn{A} ($n\in\Set{N}$) accordingly. For example, \tenfactors{\matn[i_1]{A},\matn[i_2]{A},\ldots,\matn[i_N]{A}} is a transpose of \tensor{Y}, where $(i_1,i_2,\ldots,i_N)$ is a permutation of $(1,2,\ldots,N)$. In other words, tensor transpose re-permutes the order of dimensions (modes) of \tensor{Y}. Different to matrix case, there are many ways (i.e. $N!-1$) to transpose an $N$th-order tensor.

{\bf Tensor Unfolding}. Consider the vectorization of (\ref{eqGeneralCPModel})
\begin{equation}
  \label{eqCPVector}
  \mat{y}=\sum_{j=1}^J\lambda_j\left[\matn[N]{a}_j\krp\matn[N-1]{a}_j\krp\cdots\krp\matn[1]{a}_j\right]+\mat{e},
\end{equation}
where \mat{y} and \mat{e} are respectively the vectorizations of \tensor{Y} and \tensor{E} in the proper order of dimensions. Thanks to the associativity of the Khatri-Rao product, in (\ref{eqCPVector}) we replace some Khatri-Rao products of $\matn{a}_j$ with successive $n$, $n\in\Set{N}$, by new vectors $\matn[k]{g}_j$ simultaneously for all $j=1,2,\ldots,J$, such that
\begin{equation}
  \label{eqCPVGroupDef}
  \matn[k]{g}_j=\bigkrp_{p=n_{k-1}+1}^{n_k}\matn[p]{a}_j\in\Real^{(\prod_{p=n_{k-1}+1}^{n_k}I_p)\times1},
\end{equation}
$k\in\Set{K}=\set{1,2,\ldots,K}$ with $1<K<N$, and the ascending sequence $\set{n_0,n_1,\ldots,n_K}$  forms a split of $\set{1, 2, \ldots,N}$ with $n_0=0$ and $n_K=N$. Substituting \eqref{eqCPVGroupDef} into  \eqref{eqCPVector}, we have
\begin{equation}
  \label{eqCPVectorGroup}
  \mat{y}=\sum_{j=1}^J\lambda_j\left[\matn[K]{g}_j\krp\matn[K-1]{g}_j\krp\cdots\krp\matn[1]{g}_j\right]+\mat{e}.
\end{equation}
Rewrite (\ref{eqCPVectorGroup}) in its tensor form of (\ref{eqGeneralCPModel}), we obtain a new  tensor
\begin{equation}
  \label{eqCPTensorUnfolding}
  \tenten[K]{Y}=\sum_{j=1}^J\lambda_j\matn[1]{g}_j\outerp\matn[2]{g}_j\outerp\cdots\outerp\matn[K]{g}_j+\tenten[K]{E}
\end{equation}
with the size of ${I^{\set{K}}_1\times I^{\set{K}}_2\times \cdots \times I^{\set{K}}_K}$, where $I^{\set{K}}_k=\prod_{p=n_{k-1}+1}^{n_k}I_p$, $k\in\Set{K}$. It can be observed that \tensor{Y} and \tenten[K]{Y} have exactly the same entries that are arranged in different orders. Particularly, from (\ref{eqCPVGroupDef}), by defining $\matn[k]{G}=\begin{bmatrix}
  \matn[k]{g}_1,  & \matn[k]{g}_2, & \ldots, &\matn[k]{g}_J
\end{bmatrix}$, $k\in\Set{K}$, we have
\begin{equation}
  \label{eqCPBkAn}
  \matn[k]{G}=\bigkrp\nolimits_{p=n_{k-1}+1}^{n_k}\matn[p]{A}.
\end{equation}
In other words, the tensor unfolding operation actually groups and replaces the original factors by their Khatri-Rao products. For example, \tenfactors{\matn[1]{A},\bigkrp_{k=2}^{N-2}\matn[k]{A},\matn[N]{A}\krp\matn[N-1]{A}} is a 3-way unfolding of \tensor{Y} with the dimensionality of $I_1\times (\prod_{k=2}^{N-2}I_k) \times ({I_{N-1}I_N})$; and $\tenfactors{\matn[1]{A},\ldots,\matn[N-2]{A},\matn[N]{A}\krp\matn[N-1]{A}}$ is an $(N-1)$-way unfolding. We call the above tensor unfolding  \emph{Mode Reduction} of tensors since unfolded tensors have less numbers of modes than the original ones. For simplicity we also use a notation like \tenten[N-1]{Y}=\tentenidx{1,2,\ldots,N-2,N\krp (N-1)} to denote the above unfolding where the last two modes are merged (for simplicity they are referred to as merged modes).

By using the tensor transpose and unfolding  operators introduced above, the factors \matn{A}, $n\in\Set{N}$, can be arbitrarily grouped by their Khatri-Rao products, thereby leading to different mode reduced tensors. Besides the number of modes is reduced, the mode reduced tensors have several important features which are of our particular interest:

  \subsubsection{Recoverability of original components }
  In \eqref{eqCPBkAn}, \matn[k]{G}, $k\in\Set{K}$, have very special Khatri-Rao product structures. It is known that the component matrices \matn{A}, $n\in\Set{N}$, can be estimated immediately and essentially uniquely\footnote{By essential uniqueness we mean that $\matn{\widehat{A}}=\matn{A}\mats{P}\mats{D}$  for any $\matn{\widehat{A}}$ satisfying that $\matn[k]{G}=\bigkrp_{p=n_{k-1}+1}^{n_k}\matn[p]{\widehat{A}}$, where \mats{P},  \mats{D} are any permutation matrix and nonsingular diagonal matrix, respectively,  $n\in\Set{N}$.} from \matn[k]{G}, $k\in\Set{K}$, by using a Khatri-Rao Product structure recovering procedure \cite{SPL-smbss}, i.e. the Khatri-Rao product projection (KRProj) which solves the optimization problems:
  \begin{equation}
    \label{eqKRPprojBk}
    \min_{\matn[p]{A}} \frob{\matn[k]{G}-\bigkrp\nolimits_{p=n_{k-1}+1}^{n_k}\matn[p]{A}}^2, \; k\in\Set{K}.
  \end{equation}
 Due to this fact we may decompose the mode reduced tensor \tenten[K]{Y} first such that $\tenten[K]{Y}=\tenfactors{\matn[1]{G},\matn[2]{G},\ldots,\matn[K]{G}}$ and finally recover the original components \matn{A}, $n\in\Set{N}$, by sequentially applying KRProj to the estimated component matrices \matn[k]{G}, $k\in\Set{K}$. This is the basic idea of this paper.

  In the KRProj procedure a given matrix \matn[k]{G} is optimally approximated by the Khatri-Rao product of a set of matrices with proper pre-specified size. The problem \eqref{eqKRPprojBk} generally can be solved very efficiently via SVD \cite{SPL-smbss}. A more comprehensive discussion on this topic will be detailed at the end of Section \ref{subsec:KRProj}. For notational simplicity, we denote this procedure by KRProj(\matn[k]{G}).

    \subsubsection{Non-decreasing Kruskal ranks of factors}
    First we introduce the Kruskal rank of a matrix:

  {\bf Definition 1.} The Kruskal rank of a matrix \mat{A}, i.e. $\krank{\mat{A}}$, is the largest value of $r$ such that any subset of $r$ columns of the matrix is linearly independent \cite{Kruskal1977,Sidiropoulos2000}.

  Obviously, $\krank{\mat{A}}\le\rank{\mat{A}}$ for any matrix $\mat{A}$. The following Lemma will be useful:

  \begin{lemma}[Lemma 3.3 in \cite{Stegeman2007540}]
    \label{lemma:krank}
    Consider matrices $\mat{A}\in\Real^{I_1\times J}$ and $\mat{B}\in\Real^{I_2\times J}$.  If $\krank{\mat{A}}\ge1$ and $\krank{\mat{B}}\ge1$, then $\krank{\mat{A}\krp\mat{B}}\ge\min(\krank{\mat{A}}+\krank{\mat{B}}-1,J)$.
  \end{lemma}

  From Lemma \ref{lemma:krank}, we have
  \begin{corollary}
  \label{coro_krank}
    Given a set of matrices $\matn[p]{A}\in\Real^{I_p\times J}$ with $\krank{\matn[p]{A}}\ge1$, $p=1,2,\ldots,P$, $P\ge2$, there holds that
\renewcommand{\labelenumi}{\roman{enumi}) }
    \begin{enumerate}[\IEEEsetlabelwidth{2)}]
      \item $\krank{\bigkrp_{p=1}^P\matn[p]{A}}\ge\min(J,\sum_{p=1}^P\krank{\matn[p]{A}}-(P-1)) $;

      \item $\krank{\bigkrp_{p=1}^P\matn[p]{A}}\ge \krank{\bigkrp_{p=p_0}^{p_1}\matn[p]{A}}$ where $1\le p_0\le p_1\le P$.
    \end{enumerate}

  \end{corollary}

 The proof can be found in Appendix \ref{app:ProofKrank}.

\begin{corollary}
\label{coro_krank2}
  In  Corollary \ref{coro_krank} we assume $\krank{\matn[p]{A}}\ge2$ and let $\mat{G}=\krank{\bigkrp_{p=1}^P\matn[p]{A}}$. If $P\ge J-1$, then  $\krank{\mat{G}}=\rank{\mat{G}}=J$.
\end{corollary}
\begin{IEEEproof}
If $P\ge J-1$, $\sum_{p=1}^P\krank{\matn[p]{A}}-(P-1)\ge 2P-P+1=P+1\ge J$, which means that $J\ge\rank{\mat{G}}\ge \krank{\mat{G}}\ge J$, from i) of Corollary \ref{coro_krank}. This ends the proof.
\end{IEEEproof}
Corollary \ref{coro_krank2} means that the Khatri-Rao product of at most $J-1$ matrices can form a rank-$J$ matrix, under mild conditions. From Corollary \ref{coro_krank} and \ref{coro_krank2}, the factor matrices of a mode reduced tensor, i.e. \matn[k]{G} in \eqref{eqCPBkAn}, are more likely to be of full column rank.

   \subsubsection{Less collinearity of components in merged modes}
In CPD the so-called {\emph{bottleneck}} problem that is caused by highly collinearity of some components in at least one mode \cite{PComon2009,SPL-smbss}  is a key issue that often impairs CPD algorithms. The collinearity of two vectors can be measured by
\begin{equation}
  \label{eq::corrDef}
  \rho(\mat{u},\mat{v})=\frac{|\mat{u}^T\mat{v}|}{\sqrt{\mat{u}^T\mat{u}}\sqrt{\mat{v}^T\mat{v}}}\le1,
\end{equation}
and $\rho(\mat{u},\mat{v})=1$ if and only if $\mat{u}=k\mat{v}\neq\matO$ where $k$ is a scalar. By high collinearity we mean that the collinearity measurements $\rho$ between some columns of one or more factor matrices are very close to 1, which leads to ill-conditioned factors. After tensor unfolding such that $\matn[k]{G}=\bigkrp_{p=n_{k-1}+1}^{n_k}\matn[p]{A}$ , it can be verified that for any two columns of \matn[k]{G}, say $\matn[k]{g}_{j_1}$ and $\matn[k]{g}_{j_2}$ with $j_1\neq j_2$, we have
\begin{equation}
  \label{eq::corrOfCol}
  \rho\left({\matn[k]{g}_{j_1},\matn[k]{g}_{j_2}}\right)
  =  \prod_{p=n_{k-1}+1}^{n_k}\rho\left({\matn[p]{a}_{j_1},\matn[p]{a}_{j_2}}\right)
  \le  \rho\left({\matn[p]{a}_{j_1},\matn[p]{a}_{j_2}}\right)
\end{equation}
for any $p=n_{k-1}+1,n_{k-1}+2,\ldots,n_k$. In other words, the columns of factors in merged modes after tensor unfolding become less collinear than the original ones. Particularly, once two columns in any one factor matrix, say \matn[p_0]{A}, are orthogonal, the corresponding columns in \matn[k]{G} are also orthogonal. This feature is quite helpful to improve the robustness of CPD algorithms especially when bottlenecks exist in the tensor to be decomposed. The basic trick here is that we merge the factor matrices \matn{A} which contain highly collinear columns with well-conditioned ones as possible, thereby leading to a new mode reduced tensor with well-conditioned factor matrices.

Later we will see that  the above features play very important role in the proposed method.

\subsection{A Special Case Of Mode Reduction: Matricization}
\label{subsec:MR2}
Ordinary matricization (unfolding) of tensors can be viewed as a special case of mode reduction where all modes but one are merged using their Khatri-Rao product to form 2nd-order tensors (i.e. matrices), as defined in \eqref{CPmatn}. In this case once \matn{A} and \matn{B} have been correctly estimated, all the other factors $\matn[p]{A}, p\neq n,$ can be estimated from the Khatri-Rao product projection procedure of \matn{B}, i.e. KRProj(\matn{B}), thanks to the special Khatri-Rao product structure of \matn{B} shown in \eqref{CPBn}. This allows us to estimate one factor \matn{A} with full column rank first, and then turn to the other factors to achieve CP decompositions of high-order tensors.

In \cite{SPL-smbss} we considered the case where we occasionally have some \apriori knowledge on the components (columns) of one factor, say \matn{A} with full column rank. We assume that the \apriori knowledge suffices the separability of $\matn{a}_j$ ($j=1,2,\ldots,J$) from their linear mixtures, which allows us to apply BSS methods to the unfolding matrix  $\tenmat{Y}$ to recover \matn{A} essentially uniquely due to the relationship of $\tenmat{Y}=\matn{A}\matn{B}{}^T$, i.e.
\begin{equation}
  \matn{\widehat{A}}=\Psi(\tenmat{Y})=\matn{A}\mats{P}\mats{D},
\end{equation}
 where $\Psi$ symbolically denotes a suitable BSS algorithm, and \mats{P}, \mats{D} are any permutation matrix and nonsingular diagonal matrix, respectively.
After that  let $\matn{B}=[\pinv{\matn{A}}\tenmat{Y}]^T$ and all the factors \matn[p]{A}, $p\neq n$, are estimated from KRProj($\matn{B}$). This has led to the CP decomposition method based on single mode BSS (CP-SMBSS) \cite{SPL-smbss}. As analyzed in \eqref{eq::corrOfCol}, this approach is able to overcome the bottleneck problem because high collinearity unlikely exists any more in \matn{B}.

Going a little further, we consider more general matricizations of tensors for $K$=2 in \eqref{eqCPVGroupDef}-\eqref{eqCPBkAn} such that \tenten[2]{Y}=\tenfactors{\matn[1]{G},\matn[2]{G}}, i.e.
\begin{equation}
  \label{eq:tenmat2GB}
  \tentenmat[2]{Y}{1}=\matn[1]{G}\matn[2]{G}{}^T, %\in\Real^{I_1I_2\cdots I_n\times I_{n+1}I_{n+2}\cdots I_N}.
\end{equation}
where
\begin{equation}
\begin{aligned}
  \label{eq:tenmat2G}
  \matn[1]{G}& =\bigkrp_{p=1}^n{\matn[p]{A}}\in\Real^{I_1I_2\cdots I_n\times J}, \\
  \matn[2]{G} & =\bigkrp_{p=n+1}^N\matn[p]{A}\in\Real^{I_{n+1}I_{n+2}\cdots I_N \times J},
\end{aligned}
\end{equation}
 and $1<n<N$. By incorporating the tensor transpose operator, actually we can arbitrarily split the factors \matn{A}, $n=1,2,\ldots,N$, into \emph{two} groups to form a different matricization of a given tensor.

The first advantage of \eqref{eq:tenmat2GB} compared with ordinary matricization defined in \eqref{CPmatn} lies in rank estimation of tensors, which is a very fundamental research topic in tensor analysis. If \matn{A} is underdetermined (i.e., $I_n<J$), it is generally impossible to infer the rank (i.e. the  value of $J$) from \tenmat{Y} without additional assumptions, because in this case the rank of \tenmat{Y} is also less than $J$. However, from Corollary \ref{coro_krank}, both \matn[1]{G} and \matn[2]{G} in \eqref{eq:tenmat2GB} are more likely to be of full column rank if both of them are Khatri-Rao products of multiple factor matrices, thereby probably leading to $\rank{\tentenmat[2]{Y}{1}}=J$. This feature may significantly improve the accuracy of rank estimation of tensors.

We may also extend the idea of CP-SMBSS to general matricization. For example, we apply BSS on \tentenmat[2]{Y}{1} to estimate \matn[1]{G} and then $\matn[2]{G}=[\pinv{\matn[1]{G}}\tentenmat[2]{Y}{1}]^T$, both are essentially unique. Finally all components are estimated from KRProj(\matn[1]{G}) and KRProj(\matn[2]{G}), respectively. This may lead to more robust and flexible version of CP-SMBSS, from Corollary \ref{coro_krank} and \eqref{eq::corrOfCol}.

%\cite{OGrady2005,BussgangProperty}

\section{CPD Based On Mode Reduction (MRCPD)}
\label{sec:MRCPD}
\subsection{The MRCPD Algorithms}
\label{subsec:MRCPDalgs}
Note that applying unconstrained matrix factorization to $\tentenmat[2]{Y}{1}$ in \eqref{eq:tenmat2GB} or \tenmat{Y} in \eqref{CPmatn} is unable to give desired results, because unconstrained matrix factorization suffers from rotational ambiguity, e.g.   $\tentenmat[2]{Y}{1}=[\matn[1]{G}\mat{U}][\matn[2]{G}{\mat{U}}]^T$ holds for any orthogonal matrix \mat{U} with proper size. In CP-SMBSS, BSS is employed to avoid this ambiguity by incorporating \apriori knowledge on components in one mode. In the following, we consider the case where such \apriori information is not available. As the CPD of a tensor is often essentially unique and free of rotational ambiguity under mild conditions, we run CPD algorithms to estimate the components of a mode reduced tensor at first, and then perform Khatri-Rao product projection on them to estimate the original components. As 3rd-order tensors are the simplest model with uniqueness guarantee under mild conditions \cite{Kruskal1977,Sidiropoulos2000} (see also Lemma \ref{lemma:Uniqueness} in section \ref{sec:Uniqueness}), we are interested in converting any $N$th-order tensor ($N>3$) into 3rd-order tensors in this paper, which leads to the CPD method based on mode reduction (MRCPD) of tensors listed in Algorithm \ref{algMrCPD}. In this section we assume that the involved CPDs are essentially unique. Detailed uniqueness analysis will be discussed in Section IV.
\begin{algorithm}
\caption{The General MRCPD Algorithm}
\label{algMrCPD}
\begin{algorithmic}[1]
 \REQUIRE \tensor{Y}, $J$, and a CP algorithm $\Psi$.
 \STATE Let \tenten{Y}=\tenfactors{\matn[1]{G},\matn[2]{{G}},\matn[3]{{G}}} be a 3-way unfolding of \tensor{Y}.
\STATE  Let$(\matn[1]{{G}},\matn[2]{{G}},\matn[3]{{G}})\from{\Psi}(\tenten{Y})$.
\STATE $\matn[n]{{A}}$ ($n\in\Set{N}$) are estimated via efficient Khatri-Rao product projection procedures KRProj($\matn[k]{{G}}$), $k=1,2,3$.
\RETURN $\matn[n]{\widehat{A}}, n=1,2,\cdots,N$.
\end{algorithmic}
\end{algorithm}

In \cite{n3cp::cis2012} we have considered a simpler case where $\matn[1]{G}=\matn[1]{A}$ and $\matn[2]{G}=\matn[2]{A}$. Compared with traditional $N$-way CPD methods, the MRCPD method only needs to estimate three factors at first, hence frequently unfolding to $N$ modes is avoided. This feature also makes MRCPD more easily to escape from local solutions. Note also that some excellent CPD algorithms are only developed for 3rd-order tensors, such as the self-weighted alternating trilinear decomposition method (SWATLD) \cite{swatld2001} , etc. The MRCPD method makes it possible to apply these methods seamlessly to the tensors whose orders are higher than 3.

From the analysis in Section \ref{subsec:MR2} and \cite{SPL-smbss} , once only one factor matrix with full column rank has been correctly estimated, all the other factors can then be essentially uniquely estimated. For this reason, in contrast to Algorithm \ref{algMrCPD}, we consider estimating only one factor first, say \matn[1]{G} which is of full column rank. For this purpose, we keep \matn[1]{G} unchanged, but significantly reduce the size of \matn[2]{G} and \matn[3]{G}:
\begin{itemize}
  \item {\bf Step 1}: Let \tenten{\tilde{Y}}=\tenfactors{\matn[1]{G},\matn[2]{\tilde{G}},\matn[3]{\tilde{G}}}, where \matn[2]{\tilde{G}} and \matn[3]{\tilde{G}} are factor matrices obtained by reducing the number of rows of \matn[2]{G} and \matn[3]{G}, respectively.
  \item {\bf Step 2}: Run 3-way CPD on \tenten{\tilde{Y}} to obtain \matn[1]{G};
  \item {\bf Step 3}: $\matn[3]{G}\krp\matn[2]{G}=\tentenmat{Y}{1}{}^T\pinv{\matn[1]{G}{}^T}$.
  \item {\bf Step 4}: Estimate \matn{A} ($n\in\Set{N}$) from  KRProj(\matn[1]{G}) and KRProj($\matn[3]{G}\krp\matn[2]{G}$), respectively.
\end{itemize}

The above way is quite efficient especially for large-scale problems because we can significantly reduce the size in two modes. Unfortunately it is not always the case that \matn[1]{G} is of full column rank. In such a case we may reduce the size of only one factor matrix, say \matn[3]{G}, which is generally with the largest size among the three factors, while remaining the other two unchanged. After \matn[1]{G} and \matn[2]{G} have been correctly estimated, we have
 \begin{equation}
   \label{eqMrCPD2}
   \matn[3]{G}=\tentenmat{Y}{3}\pinv{(\matn[2]{G}\krp\matn[1]{G})^T}.
 \end{equation}
In the next step, all the factor matrices \matn{A} ($n\in\Set{N}$) can be recovered from KRProj(\matn[k]{G}), $k=1,2,3$.
Note that the full column rank of $\matn[2]{G}\krp\matn[1]{G}$ is a \emph{necessary} condition of uniqueness for 3-way CPD \cite{Stegeman2007540}. Consequently, if the corresponding 3-way CPD is essentially unique, $\matn[2]{G}\krp\matn[1]{G}$ is always of full column rank. Hence we can always reduce the size in at least one mode in practice to achieve higher efficiency.

Now we discuss how to reduce the size of one mode of \tenten{Y}. Suppose that we want to reduce the size of \matn[3]{G} without changing the other factors, we consider its mode-3 matricization
\begin{equation}
  \label{eqY3mat3}
  \tentenmat{Y}{3}=\matn[3]{G}(\matn[2]{G}\krp\matn[1]{G})^T,
\end{equation}
from which we observe that reducing the rows of \matn[3]{G} is equivalent to reducing the rows of \tentenmat{Y}{3}. Hence, the following dimensionality reduction techniques may be employed:
\begin{enumerate}
\item{PCA (Truncated SVD).} Consider the truncated SVD of $\tentenmat{Y}{3}$ such that $\tentenmat{Y}{3}=\mat{UDV}^T$, where $\mat{D}\in\Real^{J\times J}$ is a diagonal matrix whose diagonal elements consist of the leading $J$ singular values of $\tentenmat{Y}{3}$. Then $\tentenmat{Y}{3}$ is updated as $\mat{V}$ by letting \matn[3]{\tilde{G}}\from$\mat{D}^{-1}\mat{U}^T\matn[3]{G}\in\Real^{J\times J}$, thereby leading to the significantly reduced size of \tenten{Y}.
\item{Fiber Sampling.} Sometimes we need to maintain the physical meaning of  original data, e.g., nonnegativity\footnote{However, it does not mean that SVD cannot be used for nonnegative data analysis. Due to the essentially uniqueness of CPD, the resulting factors \matn[1]{A} and \matn[2]{A} are essentially unique and hence can be nonnegative after adjusting the signs of their columns accordingly, no matter whether $\tenten{Y}$ and \matn[3]{G} are negative or not. See also \cite{TSP-lraNMF} for related discussion.}. In this case we can achieve dimensionality reduction by sampling the rows of  matrix $\tentenmat{Y}{3}$ \cite{CURTensorCaiafaC2010,CUR2009}, which is just equivalent to sampling the rows of \matn[3]{G}.
\end{enumerate}

Moreover, we may use the high-order SVD (HOSVD) or multilinear PCA (MPCA) methods to perform dimensionality reduction and data compression as the pre-processing step \cite{HOSVD2000,TNN2008MPCA}. It is worth noticing that the above techniques also provide an efficient CPD method for 3rd-order tensors incorporating dimensionality reduction techniques. Very often one mode, say $\matn[3]{G}=\matn[3]{A}$, can be of extremely large size. We can reduce the size of \matn[3]{G} first and then estimate \matn[3]{A} from (\ref{eqMrCPD2}). This way provides a trade-off between accuracy and efficiency and it is quite useful for large scale data.

\subsection{Khatri-Rao Product Projection (KRProj)}
\label{subsec:KRProj}
A general optimization problem to perform Khatri-Rao product projection can be formulated as
\begin{equation}
  \label{eqpikrpproj}
  \min_{\matn[k]{A}, k\in\Set{K}} \quad \frob{\mat{H}-\matn[K]{A}\krp\matn[K-1]{A}\krp\cdots\krp\matn[1]{A}}^2,
\end{equation}
where $\Set{K}=\{1,2,\ldots,K\}$. In (\ref{eqpikrpproj}), a given data matrix \mat{H} is approximated by the Khatri-Rao product of a set of matrices with specified size. In the proposed MRCPD method, $\mat{H}$ denotes the factor matrices \matn[k]{G} estimated by applying any CPD method on the mode reduced tensor \tenten[K]{Y}. By solving $K$ problems like \eqref{eqpikrpproj} (i.e. \eqref{eqKRPprojBk}) sequentially for $k\in\Set{K}$ we can estimate all factor matrices \matn{A} of the original tensor \tensor{Y}.

In (\ref{eqpikrpproj}), the columns of \matn[k]{A} ($k\in\Set{K}$) can be estimated sequentially by solving $J$ least squares problems
\begin{equation}
  \label{eqkrpj}
  \min_{\matn[k]{a}_j, k\in\Set{K}} \quad \frob{\mat{h}_j-\matn[K]{a}_j\krp\matn[K-1]{a}_j\krp\cdots\krp\matn[1]{a}_j}^2,
\end{equation}
where $j=1,2,\ldots,J$, $\mat{h}_j$ and $\matn[k]{a}_j$ are the $j$th columns of \mat{H} and \matn[k]{A}, respectively. The solutions of (\ref{eqpikrpproj}) are generally unique and can be solved, for example, by the procedure described in \cite{SPL-smbss}. On the other hand, from (\ref{eqkrpj}), we can reshape $\mat{h}_j$ such that $\mat{H}^{(j)}\approx\matn[K]{a}_j(\matn[K-1]{a}_j\krp\cdots\krp\matn[1]{a}_j)^T$, that is, $\mat{H}^{(j)}$ can be considered as the mode-$K$ unfolding of a rank-1 tensor $\tensor{H}^{(j)}$, and (\ref{eqkrpj}) is equivalent to
\begin{equation}
  \label{eqkrpTensorH}
  \min_{\matn[k]{a}_{j}, k\in\Set{K}} \quad \frob{\tensor{H}^{(j)}-\matn[1]{a}_j\outerp\matn[2]{a}_j\outerp\cdots\outerp\matn[K]{a}_j}^2.
\end{equation}
In other words, the optimal $\matn[k]{a}_j$, $k\in\Set{K}$, can be obtained by seeking the optimal rank-1 CPD of  $\tensor{H}^{(j)}$. Although this can be done by applying any standard CPD method, below we consider two relatively simple yet efficient implementations.

\subsubsection{Parallel Extraction}
Consider the mode-$n$ unfolding of $\tensor{H}^{(j)}$ and (\ref{eqkrpTensorH}) is equivalent to
\begin{equation}
  \label{eqkrpunfolding}
  \min_{\matn[k]{a}_j,\; k\in\Set{K}}\quad \frob{\tenmat[k]{H}^{(j)}-\matn[k]{a}_j\mat{v}^T}^2,
\end{equation}
where $\mat{v}=\bigkrp_{p\neq k}\matn[p]{a}_j$. Hence the optimal $\matn[k]{a}_j$ is just the left singular vector associated with the largest singular value of $\tenmat[k]{H}^{(j)}$. Consequently, all the columns of \matn[k]{A}, $k\in\Set{K}$, can be estimated uniquely and parallelly by running truncated SVD on $\tenmat[k]{H}^{(j)}$. This way may considerably benefit from parallel computations.

In order to impose specific constraints on the components we may employ the power iterations \cite{MatrixComputations} to solve (\ref{eqkrpunfolding})
\begin{equation}
  \label{eqConsPowerIter}
  \begin{aligned}
    &\matn[k]{a}_j \from \lfrac{\tenmat[k]{H}^{(j)}\mat{v}}{\frob[2]{\mat{v}}^2}, \quad
    &\mat{v} \from \lfrac{\tenmat[k]{H}^{(j)T}\matn[k]{a}_j}{\frob[2]{\matn[k]{a}_j}^2},
  \end{aligned}
\end{equation}
followed by a projection operation $\mathcal{P}$ respectively
\begin{equation}
\label{eqkrpProj}
  \matn[k]{a}_j\from \mathcal{P}(\matn[k]{a}_j), \quad \mat{v}\from\mathcal{P}(\mat{v}).
\end{equation}
For example, for nonnegative constraints $\mathcal{P}$ is element-wisely defined as
\begin{equation}
  \label{eqkrpProjNN}
  \mathcal{P}_+(x)=\max(x,0),
\end{equation}
and for sparsity constraints
\begin{equation}
  \label{eqkrpProjS}
  \mathcal{P}_\text{S}(x)=\text{sign}(x)(|x|-\lambda),
\end{equation}
where $\lambda$ is a nonnegative parameter. We repeat (\ref{eqConsPowerIter})  and (\ref{eqkrpProj}) alternatively till convergence.

\subsubsection{Tensorial Power Iterations}
Analogy to the power iteration method in the matrix case, tensorial power iterations can be derived straightforwardly by using tensor operations. From (\ref{eqkrpTensorH}) we have
\begin{equation}
  \label{eqHjXv}
  \begin{split}
  %S\tenmat[n-2]{\hat{H}}^j\mat{v}=\tensor{H}^j\times_{p-2}\matn[p]{a}_j, \text{for all}\; p\neq n.
  \tenmat[k]{{H}}^{(j)}\mat{v}= & \tensor{H}^{(j)}\times_1\matn[1]{a}_j{}^T\times_2\matn[2]{a}_j{}^T\times_{k-1}\matn[k-1]{a}_j{}^T  \\
   & \times_{k+1}\matn[k+1]{a}_j{}^T\times\cdots\times_{K}\matn[K]{a}_j {}^T\\
   \defeq & \tensor{H}^{(j)}\times_{\bar{k}}\matn[\bar{k}]{a}_j{}^T.
  \end{split}
\end{equation}
Note that $\mat{v}^T\mat{v}=\prod_{p\neq k}\matn[p]{a}_j{}^T\matn[p]{a}_j$. Then the general tensorial power iteration is
\begin{equation}
  \label{eqTPI}
  \begin{aligned}
   \matn[k]{a}_j\from \frac{\tensor{H}^{(j)}\times_{\bar{k}}\matn[\bar{k}]{a}_j{}^T}{\prod_{p\neq k}\frob[2]{\matn[p]{a}_j}^2}, \quad \matn[k]{a}_j\from\mathcal{P}(\matn[k]{a}_j),
  \end{aligned}
\end{equation}
where $\mathcal{P}$ can be (\ref{eqkrpProjNN}) or (\ref{eqkrpProjS}) to impose desired constraints. We repeat (\ref{eqTPI}) alternatively for $k=1,2,\cdots,K$ till we achieve convergence. This is actually an optimal rank-1 CPD of tensors without involving matrix inverse operations and is the extension of the power iteration method \cite{MatrixComputations} in tensor scenarios.

By repeating the above procedure for $j=1,2,\ldots,J$ all the columns of \matn[k]{A} ($k\in\Set{K}$) can be obtained, which realizes the Khatri-Rao product projection of \mat{H}.

\section{Issue of Uniqueness}
\label{sec:Uniqueness}
%From the unfolding procedure, it can be seen that the higher-order CP decompositions are always essentially unique as long as the corresponding 3-way CP decompositions are unique\footnote{To see this, just note that the corresponding Khatri-Rao products of \matn{A} form a solution to the 3rd-order tensor.}. The key point here is whether the 3-way CPD is also essentially unique if the original $N$-way CPD is. If it is not the truth, the MRCPD method may lead to very poor performance even if the original $N$-way CPD is unique. In the following this important issue will be investigated.

We use tensors in order to exploit their multi-way nature as much as possible. By mode reduction, however, $N$th-order tensors are converted into lower-order ones. The first key problem is whether mode reduction destroys the algebraic structure of  original $N$-way data and hence leads to loss of information. Or equivalently, whether MRCPD is able to give consistent results with \emph{direct CPD methods} which have no mode reduction. This is clarified by the following proposition:
\begin{theorem}
\label{th::consist}
Let \tenten[K]{Y} be a $K$th-order tensor unfolding of an $N$th-order tensor \tensor{Y} with $3\le K<N$. If both \tenten[K]{Y} and \tensor{Y} have essentially unique CPD, the MRCPD method and direct  methods theoretically give essentially the same components.
\end{theorem}

The proof of Proposition \ref{th::consist} is straightforward from the fact that the corresponding Khatri-Rao products of \matn{A} always form a solution of \tenten[K]{Y}, from the definition of the unfolding operation and the relationship \eqref{eqCPBkAn}. Obviously, if the CPD of the original $N$th-order tensor is not essentially unique, the corresponding CPD of its any mode reduced one is not either because the factors of the $N$th-order tensor are always able to form the factors of its any mode reduced tensor, thereby leading to the following proposition:

\begin{theorem}
  \label{th::K2NSufficient}
     Let \tenten[K]{Y} be any tensor unfolding of an $N$th-order tensor \tensor{Y} with $3\le K<N$. If the CPD of \tenten[K]{Y} is essentially unique, then that of \tensor{Y} is also essentially unique.
\end{theorem}
%\begin{IEEEproof}
%  Suppose that the CPD of \tensor{Y} is not unique. Then all those different decompositions form solutions for the lower-order unfolded tensor \tenten[K]{Y}, which contradicts the uniqueness of CPD of \tenten[K]{Y}.
%\end{IEEEproof}

In other words, the essential uniqueness of CPD of \tenten[K]{Y} is a sufficient condition for that of \tensor{Y}. This may simplify the exploration of new uniqueness conditions for CPD and suggest that the uniqueness conditions of 3rd-order CPD may play very important role in the uniqueness analysis.

From Proposition \ref{th::consist} and \ref{th::K2NSufficient}, the key point is whether the lower-order CPD is also essentially unique if the original $N$th-order CPD is. If it is not true, the MRCPD method may lead to very poor performance even if the original $N$th-order CPD is unique. In the following this important issue will be investigated based on the well-known uniqueness condition given by Kruskal in 1977 \cite{Kruskal1977} for 3-way tensors and then extended for $N$th-order tensors by Sidiropoulos and Bro in 2000 \cite{Sidiropoulos2000}:

\begin{lemma}[Uniqueness condition for CPD \cite{Sidiropoulos2000}]
\label{lemma:Uniqueness}
For an $N$th-order tensor \tensor{Y}=\compactcp{A}, if
\begin{equation}
\label{Kruskalcondition}
  \sum\nolimits_{n=1}^{N}\text{kr}_{\matn{A}}\ge 2J+(N-1),
\end{equation}
then the decomposition is essentially unique, where $\krank{\matn{A}}$ is the Kruskal rank of \matn{A}, $n\in\Set{N}$. (For Simplicity, we call this condition  the KSB uniqueness condition hereafter.)
\end{lemma}

Before moving on,  first we assume that $\text{kr}(\matn{A})\ge2$, $\forall{n}\in\Set{N}$, as it is a necessary condition for uniqueness of CPD \cite{Stegeman2007540}. Moreover, without loss of generality, for a given tensor \tensor{Y}=\compactcp{A} we hereafter assume that $\krank{\matn[1]{A}}\ge\krank{\matn[2]{A}}\ge\cdots\ge\krank{\matn[N]{A}}$ (Otherwise we transpose the tensor by changing the roles of \matn{A} till this condition is satisfied).

Obviously, higher-order tensors are more likely to have unique CPD than lower-order ones because the left hand side of  (\ref{Kruskalcondition}) generally increases faster than the right hand side when $N$ increases. As a result, even the original $N$th-order tensor satisfies the uniqueness condition (\ref{Kruskalcondition}), it is still possible that the corresponding lower-order CPD is not essentially unique. This problem is clarified below.

\begin{theorem}
\label{th::Kunique}
Given an $N$th-order ($N\ge4$) tensor \tensor{Y}=\compactcp{A}, if the KSB uniqueness condition is satisfied for \tensor{Y}, then it is also satisfied for the $(N-1)$th-order unfolded tensor $\tenten[N-1]{Y}=\tenfactors{\matn[1]{A},\ldots,\matn[N-2]{A},\matn[N-1]{A}\krp\matn[N]{A}}$. Consequently the CPD of \tenten[N-1]{Y} is also essentially unique.
\end{theorem}
(The proof can be found in Appendix \ref{app::ProofUNI}.)

\begin{corollary}
\label{coro:K3unfolding}
  Under the assumptions of Proposition \ref{th::Kunique}, for any $3\le K<N$ there exists at least one $K$th-order unfolded tensor \tenten[K]{Y} such that its CPD is essentially unique.
\end{corollary}

Corollary \ref{coro:K3unfolding} is obvious from Proposition \ref{th::Kunique} by letting $K=N-1, N-2, \ldots, 3$, sequentially. 
In this paper we only consider $K=3$. It can be seen that Proposition \ref{th::Kunique} also provides a way to unfold a given $N$th-order tensor to a 3rd-order one where the uniqueness is maintained. For example, if $N=5$, $J=18$, and assuming that the corresponding Kruskal ranks of the factors are 10, 9, \ldots, 6, respectively, the unfolded tensor \tenten[3]{Y}=\tenfactors{\matn[1]{A},\matn[2]{A}\krp\matn[3]{A}, \matn[4]{A}\krp\matn[5]{A}} will have essentially unique CPD. In other words, the Kruskal rank of  new factor matrices after mode reduction should be mostly balanced as possible. This result shows that the method proposed in \cite{n3cp::cis2012} may fail if all the factors are ill-conditioned. In practice, the true Kruskal ranks of factor matrices are unknown. We may use the mode ranks\footnote{The mode-$n$ rank of a tensor is the rank of its mode-$n$ matricization.} to estimate the optimal way of tensor unfoldings. Note also that, besides the KSB uniqueness condition, in \cite{UniCP3:LDL2012} some relaxed uniqueness conditions were proposed for 3-way CPD. These results simplify the validation of uniqueness.

Based on the above analysis, we may find another interesting feature of the MRCPD method. For an $N$th-order tensor, in practice we may only be interested in one factor matrix which plays a role in the $N$th-order CPD and will be used for further data analysis tasks such as clustering, classification, etc. In this case we may consider a 3rd-order CPD first and then only the factor of interest is extracted by using the Khatri-Rao product projection procedure. In this way although we actually do not perform a full CPD of the original $N$th-order tensor, the extracted factor is simply consistent with the one obtained from full $N$th-order CPD, according to Proposition \ref{th::consist}. In summary, under the KSB uniqueness condition, although we reduce the number of modes of a tensors in MRCPD, we do not lose any structural information compared with the original high-order tensor.

\section{Error bounds For Truncated CPD}
\label{sec:ErrorB}
In the above we showed that theoretically the MRCPD method is able to give the consistent results with direct methods (i.e. without mode reduction) if the corresponding CPD is essentially unique. It is worth noticing that these results are based on exact decomposition. However, in order to filter out noise and/or achieve data compression, it is common that approximate CPD is desired rather than exact one in practice. In this section we investigate the performance of the MRCPD method in truncated CPD (tCPD), which corresponds to the case where $\frob{\tensor{E}}^2> 0$ in \eqref{eqGeneralCPModel} and the best rank-$J$ approximation of \tensor{Y} is pursued, provided that the rank of \tensor{Y} is larger than $J$. Generally we have

\begin{theorem}
\label{th::KErrorBound}
Given an $N$th-order ($N\ge4$) tensor \tensor{Y} and $J$, we assume that:
\begin{enumerate}
  \item \tenfactors{\matn[1]{A},\matn[2]{A},\ldots,\matn[N]{A}} is the optimal rank-$J$ tCPD of \tensor{Y} with
\begin{equation}
  \label{eq::optCPD}
  \frob{\tensor{Y}-\tenfactors{\matn[1]{A},\matn[2]{A},\ldots,\matn[N]{A}}}=\epsilon^*,
\end{equation}
where $\matn{A}\in\Real^{I_n\times J}$;
\item \tenfactors{\matn[1]{\widetilde{A}},\matn[2]{\widetilde{A}},\ldots,\matn[N]{\widetilde{A}}} is the tCPD of \tensor{Y} obtained by using the MRCPD method such that
    \begin{enumerate}
        \item \tenfactors{\matn[1]{\widetilde{A}},\matn[2]{\widetilde{A}},\ldots,\matn[N-2]{\widetilde{A}},\mat{G}} is the optimal rank-$J$ tCPD of the unfolded tensor \tenten[N-1]{Y}=\tentenidx{1,2,\ldots,N-2,N\krp (N-1)}, where $\mat{G}\in\Real^{I_NI_{N-1}\times J}$ and $\matn{\widetilde{A}}\in\Real^{I_n\times J}$, $n<N-1$.
        \item $(\matn[N-1]{\widetilde{A}}, \matn[N]{\widetilde{A}})$ is the optimal solution of
        \begin{equation*}
          \min_{\mats[1]{C},\mats[2]{C}} \; \frob{\mat{G}-\mats[1]{C}\krp\mats[2]{C}},
        \end{equation*}
        where $\mats[1]{C}\in\Real^{I_{N}\times J},\; \mats[2]{C}\in\Real^{I_{N-1}\times J}$. Moreover, $\mat{G}-\matn[N-1]{\widetilde{A}}\krp\matn[N]{\widetilde{A}}=\mat{E}$ and $\frob{\mat{E}}=\epsilon_K$.
    \end{enumerate}
\end{enumerate}

Then we have
\begin{equation}
\begin{split}
\label{eq::ErrorBound}
  \epsilon^*\le & \frob{\tensor{Y}-\tenfactors{\matn[1]{\widetilde{A}},\matn[2]{\widetilde{A}},\ldots,\matn[N]{\widetilde{A}}}}\\
  \le & \epsilon^*+\sqrt{J}\epsilon_K.
\end{split}
\end{equation}
\end{theorem}
The proof is presented in Appendix \ref{app::proofEB}.

{\emph{Remark: In Proposition \ref{th::KErrorBound} we have assumed that the involved optimal rank-$J$ tCPD always exist. Moreover, the columns of all  factor matrices except the last one (i.e. with index $N$) are normalized with unit $L_2$ norm, as mentioned in Section \ref{subsec:CPD}.}}

Consequently, although the MRCPD method is able to give essentially the same results as direct CPD methods in perfect decomposition cases, it may achieve worse accuracy if truncated CP decompositions are performed, from Proposition \ref{th::KErrorBound}. Fortunately, the upper bound of loss of accuracy can be estimated from the Khatri-Rao product projection error approximately. If the corresponding Khatri-Rao product projection error is sufficiently small, we may conclude that the MRCPD method actually gives almost the same results with direct methods. In fact, by letting $\epsilon_K=0$  in Proposition \ref{th::KErrorBound} we have:

\begin{corollary}
  Under the assumptions of Proposition \ref{th::KErrorBound}, if there holds that $\mat{G}=\mats[N]{\widetilde{A}}\krp\matn[N-1]{\widetilde{A}}$, then the MRCPD method also gives the optimal tCPD.
\end{corollary}

 This feature allows us to determine whether we can accept a decomposition result obtained by  MRCPD or not. If the result is not acceptable, we may use the result as an initialization for the next run with different configurations (e.g., choosing a different tensor unfolding and/or CPD algorithm to perform 3rd-order CPD, etc).

\section{Simulations}
\label{sec:Simulations}
Two performance indices (PI) were used to evaluate the performance of the proposed method. The first one is the mean signal-to-interference ratio (mSIR), which is defined by
\begin{equation*}
  \label{SIR}
  \text{mSIR}(\matn{A},\matn{\widehat{A}})=\frac{1}{J}\sum_{j=1}^J 10\log_{10}{\frac{\frob[2]{\matn{a}_j}^2}{\frob[2]{\matn{a}_j-\matn{\widehat{a}}_j}^2}},
\end{equation*}
where $\matn{a}_j$, $\matn{\widehat{a}}_j$ are normalized with zero mean and unit variance, and \matn{\widehat{A}} is an estimate of $\matn{A}$ with proper re-permutation of columns $\matn{\widehat{a}}_j$. The value of mSIR  reflects how well the estimated components (source) match the true original ones.
The second PI measures the fit of  the estimated tensor to the original tensor which is defined as
\begin{equation}
  \label{fit}
  \text{Fit}(\tensor{Y},\tensor{\widehat Y})=1-\frac{\|\tensor{Y}-\tensor{\widehat Y}\|_F}{\|\tensor{Y}\|_F},
\end{equation}
where $\tensor{\widehat Y}$ is an estimate of \tensor{Y}, $\text{Fit}(\tensor{Y},\tensor{\widehat Y})=1$ if and only if $\tensor{\widehat Y}=\tensor{Y}$. For synthetic data, \tensor{Y} is the original noiseless data in order to evaluate how robust of the proposed method with respect to additive noise. All the experiments were done in MATLAB 2008a on a computer with Intel i7 3.33GHz CPU and 24GB memory running Windows 7.

\begin{figure}[!t]
    \centerline{
    \includegraphics[width=0.46\textwidth]{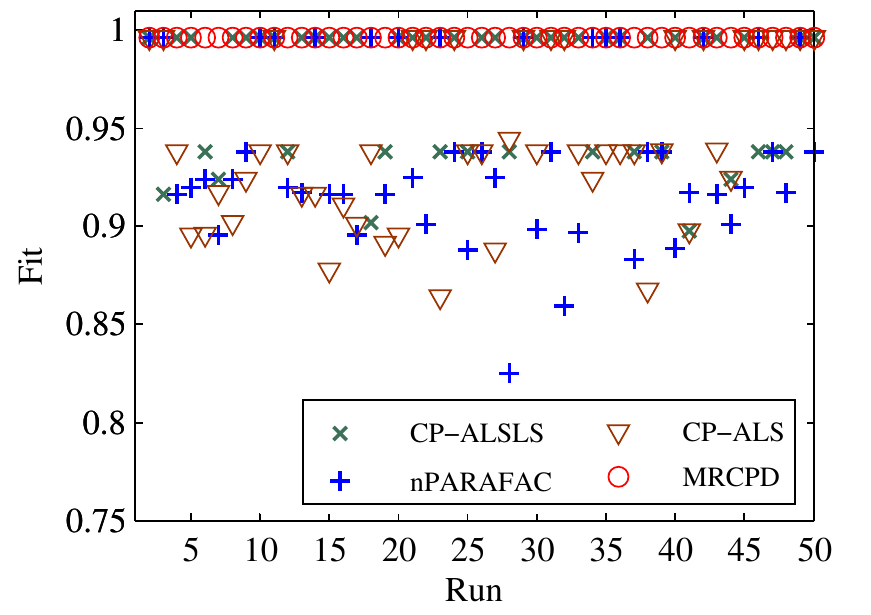}
    }
    \caption{Fit of each algorithm over 50 Monte Carlo runs. The MRCPD method consistently escaped from local solutions easily.}
    \label{figkt5randn}
  \end{figure}

{\bf Simulation 1}:
Investigation of improved efficiency and global convergence feature of the proposed method. We generated a 5th-order tensor $\tensor{Y}$ using the CP model. The elements of each factor matrix $\matn{A}\in\Real^{20\times48}$ were drawn from independent standard normal distributions, which means that each factor is underdetermined and $\sum_n\rank{\matn{A}}\le100=2J+(N-1)$. This setting makes the problem rather difficult as it is at most on the boundary of the KSB uniqueness condition. Finally independent Gaussian noise with SNR=20dB was added to the observation tensor. The proposed method was compared with the standard CP method based on ALS iterations (CP-ALS) in \cite{KoldaTensorToolbox},  the CP-ALS combined with line search (CP-ALSLS) \cite{Nion2008ALSLS}, and the PARAFAC algorithm included in the $N$-way tensor toolbox for MATLAB \cite{nwaytoolbox} (nPARAFAC, ver. 3.20). For these methods the maximum iteration number was set to 100. In the MRCPD method, we constructed a 3rd-order tensor $\tenten[3]{Y}=\tentenidx{1,3\krp 2,5\krp 4}$, and performed PCA on the mode-3 of \tenten[3]{Y}. Then we used the nPARAFAC method to perform 3-way CPD. Finally all the factors were recovered by using KRProj based on tensorial power iterations. Their performance over 50 Monte Carlo runs was detailed in TABLE \ref{tabSyn5}, where the Global Convergence Rate (GCR) evaluates the ability of escaping from local solutions  and is defined as
\begin{equation*}
  \label{eq:GCR}
  \text{GCR}=\frac{\text{Number of global convergence}}{\text{Number of Runs}}\times 100\%.
\end{equation*}
 Their corresponding Fit values are plotted in \figurename \ref{figkt5randn}, which shows that only MRCPD found the optimal Fit in all runs. From the simulation results, the MRCPD method is much more efficient and is able to escape from local minima more easily than the other methods in comparison.

\begin{figure}[!t]
    \centerline{
    \includegraphics[width=0.48\textwidth]{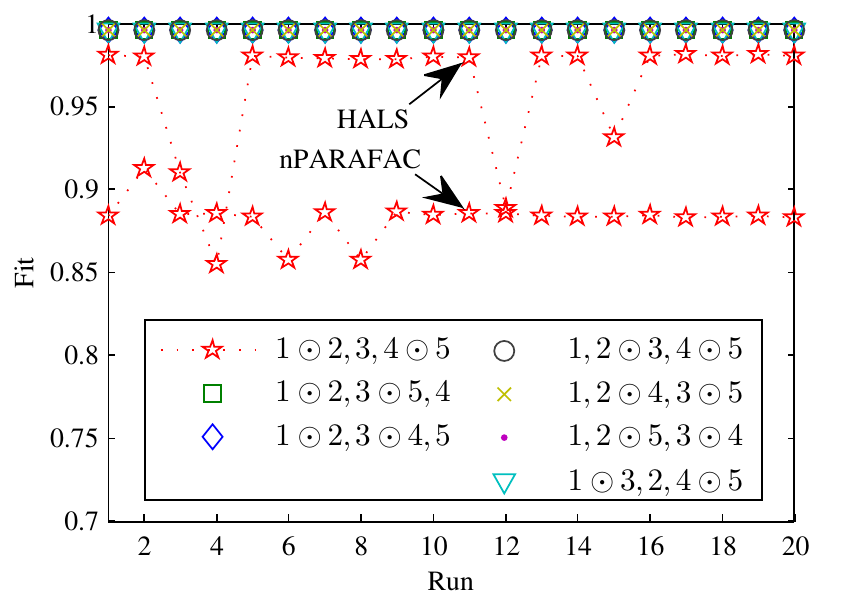}
    }
    \caption{Illustration of how tensor unfolding affects the performance in 20 Monte Carlo runs. In this experiment the nPARAFAC algorithm was used to perform 3-way CPD, except for the dot line which was obtained by using the HALS to perform 3-way CPD. It can be seen that, 1) For a given 3-way CPD algorithm, the final Fits can be different when different tensor unfoldings are used, although in the most cases there was no big difference; 2) For the same tensor unfoldings, different 3-way CPD methods lead to different final Fits, due to their own advantages, bias, and limitations.}
    \label{figkt5randnGroup}
  \end{figure}

Theoretically, the MRCPD method is able to give essentially the same results as direct CPD methods on conditon that the mode reduced tensor has unique decomposition, as analyzed in Section \ref{sec:Uniqueness}. Hence the way of tensor unfolding is important to the proposed method. To investigate the influence of different tensor unfoldings, in the following experiment we ran MRCPD under the same configuration but with 7 different randomly selected tensor unfoldings (See \figurename \ref{figkt5randnGroup} for the unfolding settings). It can be seen that among them the performance was significantly worse than the others when we used \tentenidx[Y]{1\krp2,3,4\krp5} to perform the 3-way CPD. For the other unfoldings the results are almost the same. For \tentenidx[Y]{1\krp2,3,4\krp5} we also used the hierarchical alternating least squares (HALS) algorithm  \cite{TensorHALS2009} (but without nonnegativity constraints) to replace the nPARAFAC algorithm to perform 3-way CPD, where better performance was achieved, as shown in \figurename \ref{figkt5randnGroup}. Based on this fact we may conclude: 1) For a pre-specified 3-way CPD algorithm, the final Fit may be different if different tensor unfoldings are used; 2) For the same tensor unfoldings, by selecting different 3-way CPD algorithms we may obtain different results. We think the reason is twofold. First, different tensor unfoldings may lead to tensors with different algebraic structures and features, typically, collinearity of columns. In the meanwhile, improper tensor unfoldings will result in deteriorated solutions. Second, so far existing CPD methods do not guarantee global convergence and each of them has some disadvantages and limitations. Consequently, in practice we should carefully select the way of unfolding and algorithms to perform 3-way CPD in order to achieve the best performance.

\begin{table}[!t]
%\renewcommand{\arraystretch}{1.3}
\caption{Performance comparison between the algorithms in terms of Runtime, Fit, mSIR, and GCR averaged over 50 Monte Carlo runs.}
\label{tabSyn5}
\centerline{
\begin{tabular}{c c c c c}
\hline \hline
Algorithm & Runtime(s) & Fit & mSIR & GCR\\
\hline
CP-ALS           & 59.2        &    0.94      &   52.5   &       34\%  \\
nPARAFAC     &  236.5      &   0.94      &   52.2    &      30\%  \\
CP-ALSLS      &  59.7        &   0.97       &  53.1    &    66\%     \\
MRCPD                & \bf 1.4         & \bf  1.00      &   \bf 54.9     & \bf  100\%  \\
\hline \hline
\end{tabular}
}
\end{table}

\begin{figure}[!t]
    \centerline{
    \includegraphics[width=0.5\textwidth]{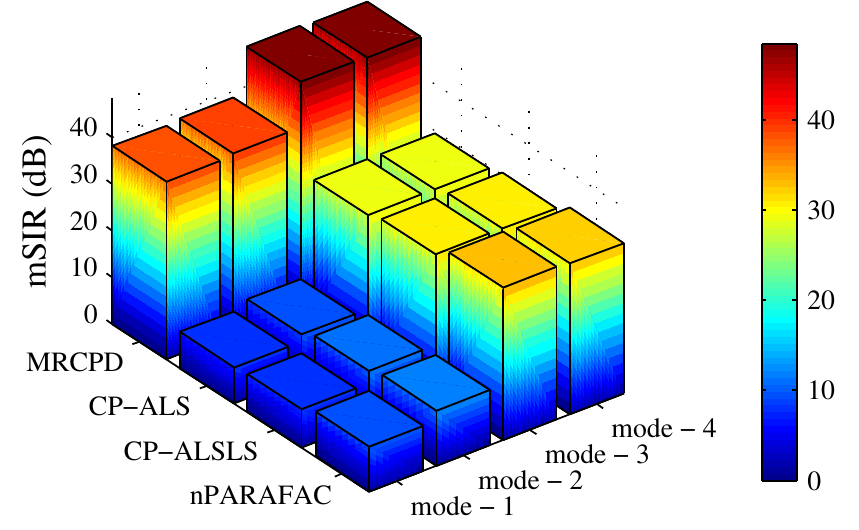}
    }
    \caption{mSIRs averaged over 50 Monte Carlo runs when they were applied to decompose a tensor with double bottlenecks. Although all the algorithms achieved very satisfying Fit, only the MRCPD method recovered all source components.}
    \label{figg12sin34sir}
  \end{figure}

{\bf Simulation 2}:
Investigation of improved uniqueness feature of the proposed method in decomposition of tensors with bottlenecks. We set $I_n=50$ ($n=1,2,\ldots,5$) and $J=5$. We used the technique adopted in \cite{PComon2009} to generate collinear components in \matn[1]{A} and \matn[2]{A}, i.e., $\matn[n]{a}_j=\mats[n,j]{v}$ for $j=1$; and $\matn[n]{a}_{j}=\matn[n]{a}_{j-1}+0.5\mats[n,j]{v}$ for $ j=2,3,\cdots,5$ and $n=1,2$, where the entries of $\mats[n,j]{v}$  were drawn from independent standard normal distributions. The components in \matn{A}, $n=3,4$, consist of 10 sine waves with slightly shifted phases but the same frequency  $f=2\text{Hz}$. That is, $\matn[3]{a}_j(t)=\sin(2\pi f t+j\frac{\pi}{50})$ and  $\matn[4]{a}_j(t)=\sin(2\pi f t+(j+5)\frac{\pi}{50})$, $j=1,2,\cdots,5$. The sampling time is from 0 to 1000 ms with the interval of 20 ms. By these configurations, the every neighboring components in \matn[n]{A} have correlations higher than 0.9. We used the same settings for all the algorithms as in Simulation 1. Their performance averaged over 50 Monte Carlo runs are detailed in TABLE \ref{tabBottleNecks}. From the table, although all the methods achieved satisfying Fit, however, only the MRCPD method recovered all the true components correctly, which shows that MRCPD is actually able to improve the uniqueness of CPD.  See \figurename \ref{figg12sin34sir} for their mSIRs averaged over 50 Monte Carlo runs.

\begin{table}[!t]
%\renewcommand{\arraystretch}{1.3}
\caption{Performance of the algorithms when they were applied to decompose a tensor with bottlenecks. Although all algorithms achieved very satisfying Fit of 1.0, only the MRCPD method recovered all the true components.}
\label{tabBottleNecks}
\centerline{
\begin{tabular}{c c c c c}
\hline \hline
Algorithm & CP-ALS & nPARAFAC & CP-ALSLS & MRCPD \\
\hline
Fit                       &1.0 & 1.0 & 1.0 & 1.0 \\
Runtime(s)          & 16.1 & 36.7 & 39.8 & \bf11.0 \\
mSIR(dB)                   &18.9 & 21.6 & 20.0 & \bf43.3 \\
\hline \hline
\end{tabular}
}
\end{table}

\begin{figure}[!t]
    \centerline{
    \subfloat[]{
    \includegraphics[width=0.5\textwidth]{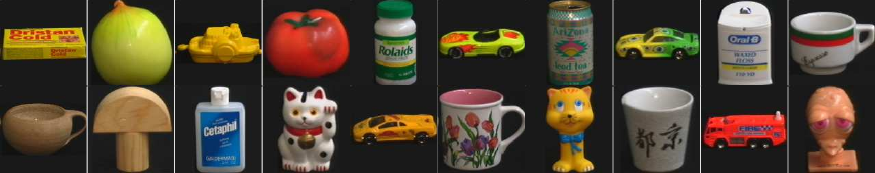}}
    }
    \centerline{
    \subfloat[]{
    \includegraphics[width=0.5\textwidth]{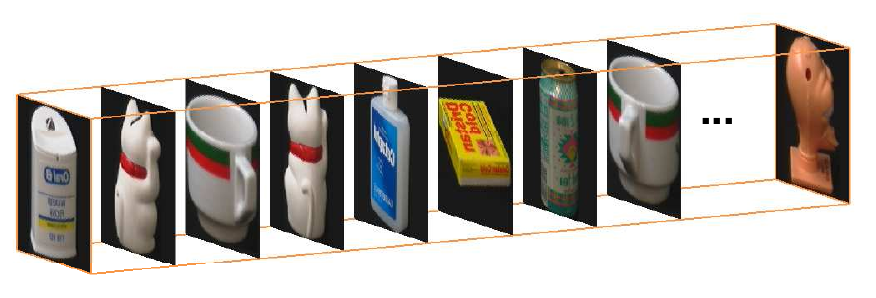}}
    }
    \caption{(a) One pose of the first 20 objects of the COIL-100 database. (b) Illustration of how to generate the observation tensor.}
    \label{figCOIL100-20}
  \end{figure}

\begin{figure}[!t]
    \centerline{
    \includegraphics[width=0.46\textwidth]{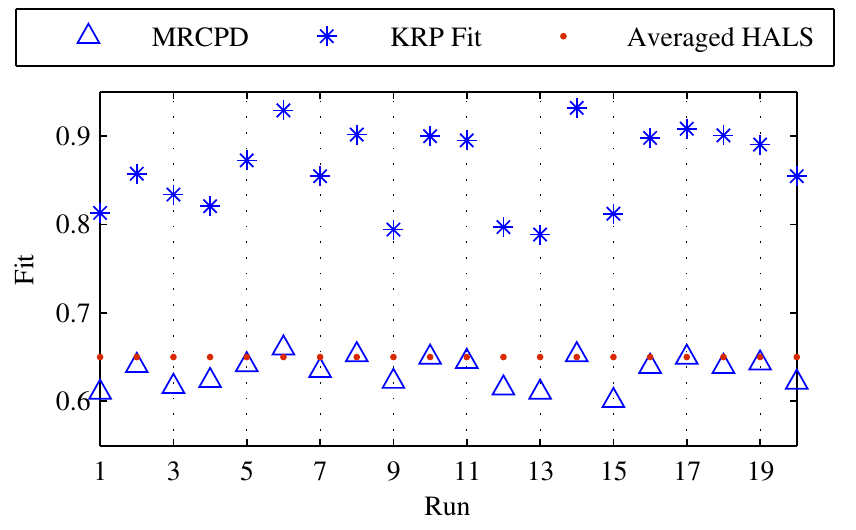}
    }
    \caption{Illustration of how the Khatri-Rao product projection accuracy (KRP Fit) affects the final Fit of the MRCPD method in the case of truncated CPD. Once the KRP Fit is satisfactory, the total Fit obtained by MRCPD is comparable with direct methods. Here KRP Fit is defined as $\epsilon_K$, see Proposition \ref{th::KErrorBound}.}
    \label{figCOILKrp}
  \end{figure}

{\bf Simulation 3}:
We applied the proposed methods to real image data analysis, namely the COIL-100 database \cite{coil100}. The COIL-100 database consists of 7200 color images of 100 objects, and 72 images per object which were taken from 72 different angles. For simplicity, we selected the first 20 objects for test, and each image was scaled with the size of $128\times128$ (See \figurename \ref{figCOIL100-20}). Then a tensor \tensor{Y} with the size of $128\times128\times3\times1440$ was generated. We set rank $J=10$ for all methods. In the MRCPD method, $\tenten[3]{Y}=\tentenidx{1;4;2\krp3}$ and we randomly sampled 100 fibers in mode-3 in each run. We used the HALS method \cite{TensorHALS2009} to perform 3-way CPD. For the other methods the maximum iteration number was set to 100. Finally, the factor \matn[4]{A} was used as features to cluster the original images. As $K$-means is prone to be influenced by initial centers of clusters, we replicated $K$-means 20 times for each method. %The widely used performance index Accuracy (\%) was adopted to evaluate the clustering results \cite{GNMF2011PAMI}.
See TABLE \ref{tabCOIL100} for their performance over 20 Monte Carlo runs. From the table, we see that the MRCPD method achieved the best clustering accuracy and it was significantly faster than the other methods. As this example is a typical truncate CPD problem, we investigated how the Khatri-Rao product projection accuracy affected the final Fit of the MRCPD method, see \figurename \ref{figCOILKrp}. It can be seen that, in truncated CPD, the Fit obtained by the MRCPD method was slightly worse than the other direct methods without mode reduction. However, once the Khatri-Rao product projection accuracy was satisfactory, the MRCPD was able to give very good results, which is consistent with Proposition \ref{th::KErrorBound}.

\begin{table}[!t]
%\renewcommand{\arraystretch}{1.3}
\caption{Performance comparison of the four CPD methods when they were applied to real image data clustering over 20 Monte Carlo runs.}
\label{tabCOIL100}
\centerline{
\begin{tabular}{c c c c }
\hline \hline
Algorithm & Fit & Runtime (s) & Accuracy \\
\hline
CP-ALS            & 0.67     & 564.2      & 63.8\%  \\
HALS                & 0.65     & 506.8   & 65.5\%    \\
nPARAFAC       & 0.67     & 746.5   & 65.4 \%  \\
MRCP               &0.64      & \bf 44.6    & \bf 67.0 \%   \\
\hline \hline
\end{tabular}
}
\end{table}
%CP-ALSLS: 0.67 & 1427.5 & 64.6

\section{Conclusion}
\label{sec:Conclusion}
Existing Canonical Polyadic decomposition (CPD) methods are often based on alternating least squares (ALS) iterations. In ALS methods, we have to unfold the observation tensor frequently, which is one major performance bottleneck of CPD. To overcome this problem, in this paper we proposed the concept of mode reduction in CPD and based on it we developed a new method to perform CPD of high-order tensors with $N\ge3$, efficiently incorporating dimensionality reduction techniques. In this way, frequently unfolding with respect to $N$ mode is avoided and the new method can escape from local solutions more easily due to the significantly reduced complexity of models. Particularly, The proposed method is promising to overcome the bottleneck problem caused by high collinearity of components, because after proper mode reduction collinearity is unlikely exist any more. Moreover, a full $N$-way CPD may be avoid but without loss of structural information of data, if only partial factors are of interest. The essential uniqueness of CPD of  mode reduced tensors was also theoretically investigated. Simulations confirmed the efficiency and validity of the proposed method.

The special algebraic structure (i.e., Khatri-Rao product structure) plays a key role in CPD of tensors, which shows that the factors in CP model have very tight intrinsic connections. In fact, once one factor with full column rank has been determined, all the other factors are also determined. Hence the classical ALS routine perhaps is not the best choice for CPD. This phenomenon may motivate us to develop more efficient and robust CPD algorithms in the future.

\appendices

\section{Proof of Lemma \ref{lemma:krank}}
\label{app:ProofKrank}
\begin{IEEEproof}
  We use mathematical induction. From Lemma \ref{lemma:krank}, we have $\krank{\matn[1]{A}\krp\matn[2]{A}}\ge\min(J,\krank{\matn[1]{A}}+\krank{\matn[2]{A}}-1)$, that is, i) holds for $P=2$.
  Now we assume that the statement i) is true for $P=K$ ($K\ge2$), that is,
  \begin{equation}
  \label{proof_induction_1}
    \krank{\bigkrp_{p=1}^K\matn[p]{A}}\ge\min(J,\sum\nolimits_{p=1}^K\krank{\matn[p]{A}}-(K-1)).
  \end{equation}
  In the following we will show that the statement is also true for $P=K+1$.
  \begin{equation*}
      \begin{aligned}
    \krank{\bigkrp_{p=1}^{K+1}\matn[p]{A}} = & \krank{\matn[K+1]{A}\krp[\matn[K]{A}\krp\matn[K-1]{A}\cdots\krp\matn[1]{A}]}\\
         \ge & \min(J,\krank{\matn[K+1]{A}}+\krank{\bigkrp\nolimits_{p=1}^{K}\matn[p]{A}}-1) \\
         \ge  & \min (J, \krank{\matn[K+1]{A}}-1 \\ & \qquad + \min(J,\sum\nolimits_{p=1}^K\krank{\matn[p]{A}}-(K-1))) \\
         \ge & \min(J,\krank{\matn[K+1]{A}}-1+J,\sum_{p=1}^{K+1}\krank{\matn[p]{A}}-K)  \\
         \ge & \min(J,\sum\nolimits_{p=1}^{K+1}\krank{\matn[p]{A}}-[(K+1)-1]),
      \end{aligned}
  \end{equation*}
  thereby showing that i) also holds for $P=K+1$. This ends the proof of statement i).
  In the above we have used that $\krank{\matn[p]{A}}\ge1,\; \forall n$, which leads to $J\le J+\krank{\matn[K+1]{A}}-1$.

  Regarding ii), note that
    \begin{equation*}
      \begin{aligned}
    \krank{\bigkrp_{p=1}^{P}\matn{A}}
         \ge & \min(J,\krank{\matn[P]{A}}+\krank{\bigkrp\nolimits_{p=1}^{P-1}{\matn[p]{A}}}-1) \\
        = & \min(J, \krank{\bigkrp\nolimits_{p=1}^{P-1}{\matn[p]{A}}}+(\krank{\matn[P]{A}}-1))  \\
        \ge & \krank{\bigkrp\nolimits_{p=1}^{P-1}{\matn[p]{A}}}.
      \end{aligned}
  \end{equation*}
  Similarly, we can prove that $\krank{\bigkrp_{p=1}^{P}\matn{A}}\ge \krank{\bigkrp\nolimits_{p=2}^{P}{\matn[p]{A}}}$. Following these procedures, we  prove ii).

The proof is complete.
\end{IEEEproof}

\section{Proof of Proposition 1}
 \label{app::ProofUNI}

%We have assume that $\krank{\matn[1]{A}}\ge\krank{\matn[2]{A}}\ge\cdots\ge\krank{\matn[N]{A}}\ge2$ as it is a necessary condition for uniqueness \cite{Stegeman2007540}. Consequently $\krank{\matn[n_1]{A}\krp\matn[n_2]{A}}\ge\min(\krank{\matn[n_1]{A}}+\krank{\matn[n_2]{A}}-1,J)$ holds for any $n_1\neq n_2$, from Lemma 2.

\begin{IEEEproof}[Proof of Proposition \ref{th::Kunique}] Consider that $\tenten[N-1]{Y}=\tenfactors{\matn[1]{A},\ldots,\matn[N-2]{A},\matn[N]{A}\krp\matn[N-1]{A}}$ where the two factors with the minimum Kruskal ranks are merged into one mode. Let $t=\krank{\matn[N]{A}\krp\matn[N-1]{A}}$. We only need to show that $\sum_{k=1}^{N-2}\krank{\matn[k]{A}}+t\ge 2J+(N-2)$, a sufficient condition of essential uniqueness of CPD of \tenten[N-1]{Y}.

From Lemma \ref{lemma:krank}, $t\ge\min(\krank{\matn[N-1]{A}}+\krank{\matn[N]{A}}-1,J)$. There are two cases:

1) $\krank{\matn[N-1]{A}}+\krank{\matn[N]{A}}-1 > J$, i.e. $t=J$ and
\begin{equation}
\label{eqAppCase1}
  2\krank{\matn[N-1]{A}}\ge\krank{\matn[N-1]{A}}+\krank{\matn[N]{A}}\ge J+2.
\end{equation}
Note that $\sum_{n=1}^{N-2}\krank{\matn[n]{A}}\ge (N-2)\krank{\matn[N-1]{A}}$. From (\ref{eqAppCase1}) and $t=J$, we have
\begin{equation}
\begin{split}
   &\sum\nolimits_{n=1}^{N-2}\krank{\matn[n]{A}}+t-[2J+(N-2)] \\
   \ge& (N-2)\frac{J+2}{2}+J-[2J+(N-2)] \\
   \ge&\frac{1}{2}(N-4)J \ge0,\\
\end{split}
\end{equation}
that is, $\sum_{k=1}^{N-2}\krank{\matn[k]{A}}+t\ge2J+(N-2)$.

2) $\krank{\matn[N-1]{A}}+\krank{\matn[N]{A}}-1 \le J$, and hence $t\ge\krank{\matn[N-1]{A}}+\krank{\matn[N]{A}}-1$. Moreover,
\begin{equation}
  \begin{split}
    &\sum\nolimits_{n=1}^{N-2}\krank{\matn{A}}+t \\
    \ge & \sum\nolimits_{n=1}^{N-2}\krank{\matn{A}}+(\krank{\matn[N-1]{A}}+\krank{\matn[N]{A}}-1) \\
    = & \sum\nolimits_{n=1}^{N}\krank{\matn{A}}-1 \\
    \ge & 2J+(N-2).
  \end{split}
\end{equation}
The proof is complete.
\end{IEEEproof}

\section{Proof of Proposition \ref{th::KErrorBound}}
\label{app::proofEB}
\begin{IEEEproof}
  The first inequality of (\ref{eq::ErrorBound}) is straightforward. We focus on the second one. Let  $\mat{H}=\matn[N-2]{\widetilde{A}}\krp\matn[N-3]{\widetilde{A}}\cdots\krp\matn[1]{\widetilde{A}}$. From $\matn{\widetilde{a}}_j{}^T\matn{\widetilde{a}}_j=1$ for any $j=1,2,\ldots,J$ and $n\neq N-1$, we have
  \begin{equation}
    \label{eq::Fnorm}
    \begin{split}
      & \frob{\mat{H}} = \sqrt{\trace{\mat{H}^T\mat{H}}}\\
                              = & \sqrt{\trace{(\matn[1]{\widetilde{A}}{}^T\matn[1]{\widetilde{A}})\hdp(\matn[2]{\widetilde{A}}{}^T\matn[2]{\widetilde{A}})\cdots\hdp(\matn[N-2]{\widetilde{A}}{}^T\matn[N-2]{\widetilde{A}})}} \\
                              = & \sqrt{J}.
    \end{split}
  \end{equation}
  Note that
  \begin{equation}
    \nonumber
    \begin{split}
           &\frob{\tensor{Y}-\tenfactors{\matn[1]{\widetilde{A}},\matn[2]{\widetilde{A}},\ldots,\matn[N]{\widetilde{A}}}}\\
      =  & \frob{\tentenmat[N-1]{Y}{N-1}-(\matn[N]{\widetilde{A}}\krp\matn[N-1]{\widetilde{A}})\mat{H}^T} \\
      =  & \frob{\tentenmat[N-1]{Y}{N-1}-(\mat{G}-\mat{E})\mat{H}^T} \\
      \le & \frob{\tentenmat[N-1]{Y}{N-1}-\mat{G}\mat{H}^T}+\frob{\mat{E}\mat{H}^T} \\
      = & \frob{\tenten[N-1]{Y}-\tenfactors{\matn[1]{\widetilde{A}},\matn[2]{\widetilde{A}},\ldots,\matn[N-2]{\widetilde{A}},\mat{G}}}+\epsilon_K\frob{\mat{H}} \\
      \le & \frob{\tenten[N-1]{Y}-\tenfactors{\matn[1]{A},\matn[2]{A},\ldots,\matn[N-2]{A},\matn[N]{A}\krp\matn[N-1]{A}}}\\
           & +\epsilon_K\sqrt{J} \\
      =  & \frob{\tensor{Y}-\compactcp{A}} +\epsilon_K\sqrt{J}\\
      \le & \epsilon^*+\sqrt{J}\epsilon_K.
    \end{split}
  \end{equation}
In the above, \tentenmat[N-1]{Y}{N-1} denotes the mode-$(N-1)$ matricization of mode reduced tensor \tenten[N-1]{Y}.

\end{IEEEproof}
%\section*{Acknowledgment}

\section*{Acknowledgment}
The authors sincerely thank the associated Editor and the anonymous reviewers for their insightful comments and suggestions that led to the present improved version of the original manuscript.

\bibliographystyle{IEEEtran}

%\bibliography{\bibpath/tensor,\bibpath/generalmath,\bibpath/all,\bibpath/optimization,\bibpath/nmf,\bibpath/database,\bibpath/zgx,\bibpath/TNN_BSS}

\begin{IEEEbiography}[{\includegraphics[width=1in,height=1.25in,clip,keepaspectratio]{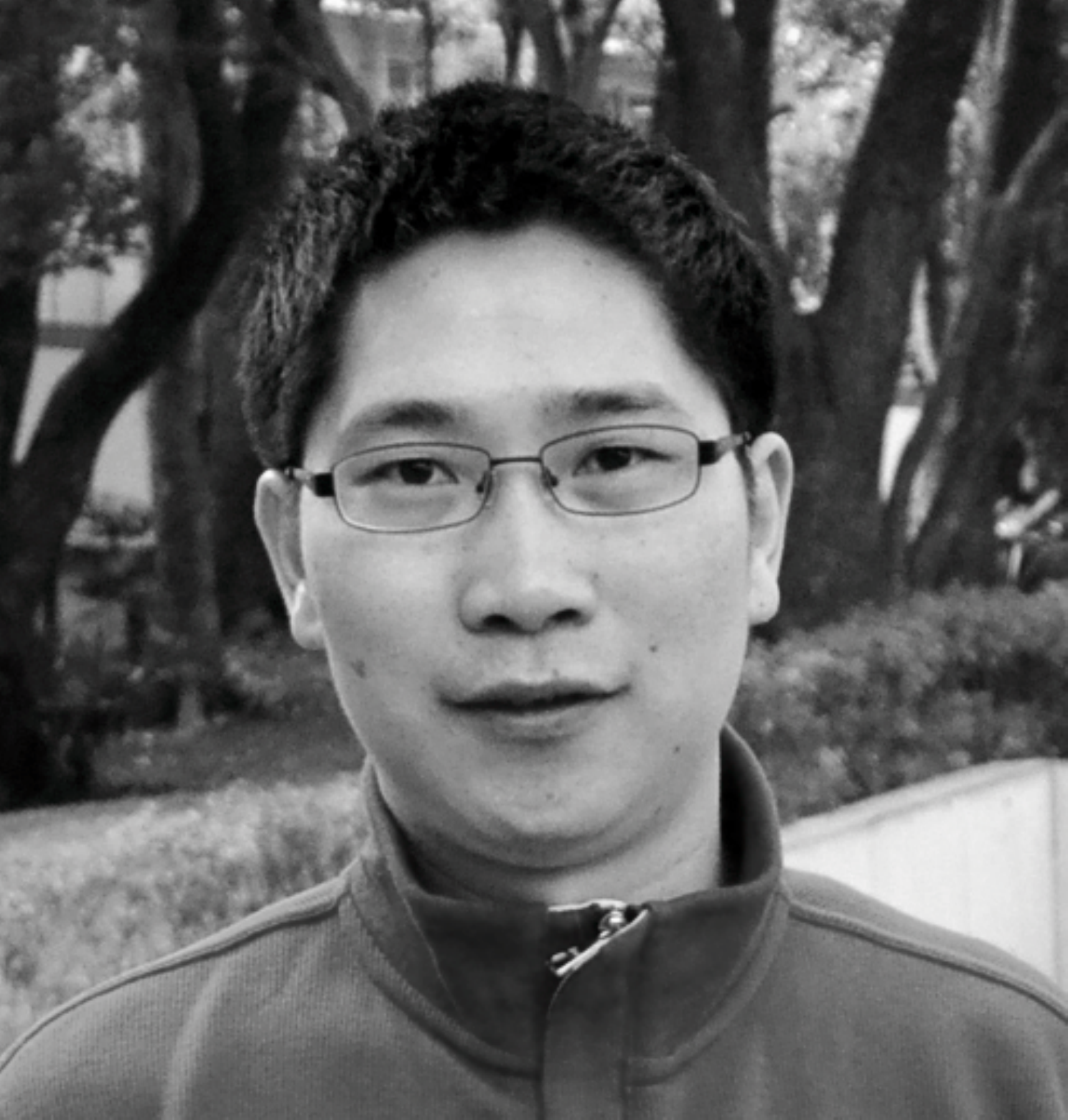}}]{Guoxu Zhou} received the Ph.D degree in intelligent signal and information processing from South China University of Technology, Guangzhou, China, in 2010. He is currently a research scientist of the laboratory for Advanced Brain Signal Processing, at RIKEN Brain Science Institute (JAPAN). His research interests include statistical signal processing, tensor analysis, intelligent information processing, and machine learning.
\end{IEEEbiography} 
\begin{IEEEbiography}[{\includegraphics[width=1in,height=1.25in,clip,keepaspectratio]{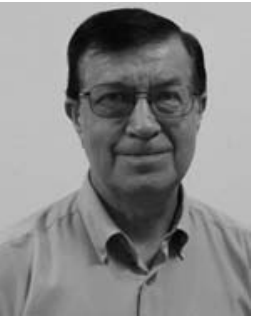}}]{Andrzej Cichocki} received the M.Sc. (with honors), Ph.D. and Dr.Sc. (Habilitation) degrees, all in electrical engineering. from Warsaw University of Technology (Poland). Since 1972, he has been with the Institute of Theory of Electrical Engineering, Measurement and Information Systems, Faculty of Electrical Engineering at the Warsaw University of Technology, where he obtain a title of a full Professor in 1995. He spent several years at University Erlangen-Nuerenberg (Germany), at the Chair of Applied and Theoretical Electrical Engineering directed by Professor Rolf Unbehauen, as an Alexander-von-Humboldt Research Fellow and Guest Professor. In 1995-1997 he was a team leader of the laboratory for Artificial Brain Systems, at Frontier Research Program RIKEN (Japan), in the Brain Information Processing Group. He is currently the head of the laboratory for Advanced Brain Signal Processing, at RIKEN Brain Science Institute (JAPAN). He is author of more than 250 technical papers and 4 monographs (two of them translated to Chinese).
\end{IEEEbiography} 
\begin{IEEEbiography}[{\includegraphics[width=1in,height=1.25in,clip,keepaspectratio]{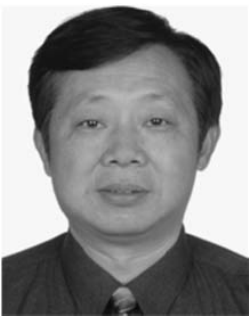}}]{Shengli Xie}(M'01-SM'02) received the M.S. degree in mathematics from Central China Normal University, Wuhan, China, in 1992, and the Ph.D. degree in control theory and applications from the South China University of Technology, Guangzhou, China, in 1997.

He is the Director of the Laboratory for Intelligent Information Processing (LIIP) and a Full Professor
with the Faculty of Automation, Guangdong University of Technology, Guangzhou. He has authored
or co-authored two monographs and more than 80 scientific papers published in journals and conference proceedings, and is a holder or joint holder of 12 patents. His current research interests include automatic control and signal processing, with a focus on blind signal processing and image processing.
\end{IEEEbiography}

\end{document}